\documentclass[12pt,a4paper]{article}
\usepackage{amsmath, amssymb}
\usepackage{amsfonts, ifthen, theorem}
\usepackage[usenames]{color}
\usepackage{mathrsfs}

\usepackage{fullpage}

\theoremheaderfont\scshape
\newtheorem{theorem}{Theorem}[section]
\newtheorem{lemma}[theorem]{Lemma}

\newtheorem{definition}[theorem]{Definition}

\theorembodyfont{\rmfamily}
\newtheorem{remark}[theorem]{Remark}
\newtheorem{remarks}[theorem]{Remarks}

\newtheorem{assumption}[theorem]{Assumption}

\renewcommand\thefootnote{\fnsymbol{footnote}}

\newcommand\R{\mathbb R}
\newcommand\N{\mathbb N}
\newcommand\F{\mathscr F}
\newcommand\E[2][\P]{\mathbb E_{#1}\left[ #2\right]}
\renewcommand\P{\mathbb P}
\newcommand\Q{\mathbb Q}

\renewcommand\d{\mathrm d}
\newcommand\qed{\hfill$\Box$}
\newcommand\ip[2]{\langle #1,#2 \rangle}
\newcommand\adm{\textup{adm}}
\newcommand\suprep{\overline\pi}
\newcommand\radnik[1][\mathbb Q]{\frac{\d #1}{\d\mathbb P}}
\newcommand\hatMV{\widehat M_V}
\newcommand\normcl[2][\Q]{\overline{#2}^{L^1(#1)}}
\newcommand\ind{1\hspace{-2.5mm}1}
\newcommand\cone{\operatorname{cone}}
\newcommand\dom{\operatorname{dom}}

\newcommand{\nohyphens}{\hyphenpenalty=10000\exhyphenpenalty=10000\relax}%
\makeatletter
\renewcommand{\@seccntformat}[1]{\csname the#1\endcsname.\hspace{1em}}%
\renewcommand\section{\@startsection{section}{1}{\z@}%
                {-3.5ex \@plus -1ex \@minus -.2ex}%
                {2.3ex \@plus.2ex}%
                {\setcounter{equation}0\bfseries\nohyphens}}%
\makeatother

\newenvironment{proof}[1][]{\noindent\textit{Proof#1.} }{\vskip\baselineskip}

\begin{document}

{\centering {\renewcommand\thefootnote{}\large\bfseries On
utility-based super-replication prices\\ of contingent claims with
unbounded payoffs\footnote{The authors gratefully acknowledge
support from EPSRC grant no. GR/S80202/01.}

}

\par
\vspace{1em} {\scshape Frank Oertel}\par\itshape Department of
Mathematics\par University College Cork\par \vspace{1em} {\scshape
Mark Owen}
\par\itshape Department of Actuarial Mathematics and
  Statistics\\ Heriot-Watt University

}

\renewcommand{\abstractname}{}
\begin{abstract}\renewcommand\thefootnote{}\footnotesize
Consider a financial market in which an agent trades with
utility-induced restrictions on wealth. For a utility function
which satisfies the condition of reasonable asymptotic elasticity
at $-\infty$ we prove that the utility-based super-replication
price of an unbounded (but sufficiently integrable) contingent
claim is equal to the supremum of its discounted expectations
under pricing measures with finite {\it loss-entropy}. For an
agent whose utility function is unbounded from above, the set of
pricing measures with finite loss-entropy can be slightly larger
than the set of pricing measures with finite entropy. Indeed, the
former set is the closure of the latter under a suitable weak
topology.

Central to our proof is the representation of a cone $C_U$ of
utility-based super-replicable contingent claims as the polar cone
to the set of finite loss-entropy pricing measures. The cone $C_U$
is defined as the closure, under a relevant weak topology, of the
cone of all (sufficiently integrable) contingent claims that can
be dominated by a zero-financed terminal wealth.

We investigate also the natural dual of this result and show that
the polar cone to $C_U$ is generated by those separating measures
with finite loss-entropy. The full two-sided polarity we achieve
between measures and contingent claims yields an economic
justification for the use of the cone $C_U$, and an open question.

\footnote{{\it AMS 2000 subject classifications.} 1B16, 46N10,
60G44} \footnote{{\it Key Words and phrases.} super-replication,
incomplete markets, contingent claims, duality theory, weak
topologies}
\end{abstract}

\section{Introduction}

Consider a financial market where the discounted prices of $d$
risky assets are modelled over a finite time interval $[0,T]$ by
an $\R^d$-valued semimartingale $S=(S_t)_{0\le t\le T}$, on a
filtered probability space $(\Omega,\F,(\F_t)_{t\in[0,T]},\P)$
satisfying the usual conditions of right continuity and
saturatedness. A portfolio on such a market can be represented by
a pair $(x,H)$, consisting of an initial wealth $x\in\R$, and a
predictable, $S$-integrable process $H$ representing the holdings
of the $d$ risky assets. It is assumed that, at any time, all
remaining wealth is invested in the numeraire. The discounted
wealth process corresponding to the portfolio $(x,H)$ is defined
by $X_t^{(x,H)}:=x+\int_0^tH_u\d S_u$.

Two important theoretical concepts within the above framework for
models of financial markets are those of ``No Arbitrage'' and
completeness. An arbitrage opportunity is defined as a trading
strategy $H$ such that $X_T^{(0,H)}\ge 0$, $\P$-a.s. and such that
$\P(X_T^{(0,H)}>0)>0$. A model is usually said to satisfy the
condition of No Arbitrage if there does not exist an {\it
admissible} trading strategy which is an arbitrage opportunity.
The condition on $H$ of admissibility is the requirement that the
wealth process $X^{(0,H)}$ is uniformly bounded below by a
constant; ruling out such processes is one way to disallow the use
of doubling strategies.

In the celebrated papers \cite{DelbScha94,DelbScha98} the
condition of No Arbitrage was weakened to that of No Free Lunch
with Vanishing Risk (NFLVR), and it was shown that a model
satisfies NFLVR if and only if the set, $M_\sigma^e$, of
equivalent $\sigma$-martingale measures is non-empty.

A model satisfying NFLVR is then said to be complete if the set
$M_\sigma^e$ is a singleton (i.e. if $M_\sigma^e=\{\Q\}$). In a
complete market model it is possible, with the use of an
admissible trading strategy, to replicate and thereby uniquely
price all contingent claims in $L^1(\Q)$ with payoffs bounded from
below (see \cite[Theorem 5.16]{DelbScha98}). In contrast to this,
if a market is incomplete there exist contingent claims with
payoffs bounded from below which cannot be replicated by
admissible trading strategies. For such contingent claims there
exists an open interval of arbitrage-free prices, rather than a
unique price.

Given a general contingent claim with payoff $X$, it is easy to
see that an upper bound for the interval of arbitrage free prices
is given by the super-replication price
\begin{equation}\label{eqn:defsuprep}
\suprep(X) := \inf\{x \in \R : \text{there is an admissible } H
\text{ such that } X \le X_T^{(x,H)}\}.
\end{equation}

As a special case of \cite[Theorem 5.5]{DelbScha98} we know that
for a contingent claim with payoff $X$ bounded from below, the
super-replication price is in fact the least upper bound for the
interval of arbitrage free prices, in other words
\begin{equation}\label{eqn:upeqlow}
\suprep(X) = \sup_{\Q\in M_\sigma^e}\E[\Q]{X}.
\end{equation}
However, for contingent claims with payoffs unbounded from below,
admissible trading strategies are unsuitable for
super-replication, and this dual representation of the
super-replication price does not always hold. Indeed,
\cite[Example 8]{BiagFrit03} constructs a market model and a
contingent claim with payoff $X$ unbounded from below such that
\begin{equation}\label{eqn:ineq}
  \suprep(X)>\sup_{\Q\in M^e_{\sigma}}\E[\Q]X.
\end{equation}
The intuitive reason for the breakdown of \eqref{eqn:upeqlow} is
that the cone
\[ K^{\adm}:=\left\{X_T^{(0,H)}:H\text{ is admissible}\right\}. \]
is not closed with respect to a weak topology induced by the set
of pricing measures.

It is useful at this point to extend slightly the definition of a
super-replication price to allow terminal wealths from an
arbitrary convex cone $K\subseteq L^0(\P)$. Let
\begin{equation*}
\suprep(X;K):=\inf\{x\in\R:X\le x+Y\textup{ for some }Y\in K\}.
\end{equation*}
Of course, $\suprep(X)=\suprep(X;K^{\adm})$. Note that if $K$ is a
solid cone in a subspace $F$ of $L^0(\P)$ (i.e. $X \in F$ and
$X\le Y\in K$ implies $X\in K$) then
\[ \suprep(X;K):=\inf\{x\in\R: X - x\in K\}. \]
If $K$ is not solid then we may of course replace $K$ by the
smallest solid cone containing $K$ without affecting the
super-replication price.

We are now able to formulate the following natural question with
\eqref{eqn:ineq} in mind: Given an arbitrary convex cone $K$ of
contingent claims, is it possible to find a minimal solid, closed
convex cone $C\supseteq K$, and a suitable set $M$ of pricing
measures, such that
\begin{equation}\label{eqn:arbcone} \suprep(X;C)=\sup_{\Q\in
M}\E[\Q]X? \end{equation}

A positive answer to this question was given in \cite{BiagFrit03}.
In this paper, preferences of an investor were incorporated in the
construction of a weakly closed, enlarged cone $C$ by means of the
convex conjugate of the investor's utility function. The set $M$
of measures consisted of those absolutely continuous separating
measures with finite entropy. A dual representation of the form
\eqref{eqn:arbcone} was obtained for utility functions which are
bounded from above. This result has since been extended in
\cite{Owen03} to unbounded utility functions with Reasonable
Asymptotic Elasticity at both $-\infty$ and $+\infty$, with
subsequent alternative proofs given in \cite{BiagFrit04},
\cite{BiagFrit05} and \cite{OertOwen05}. In this article we show
that the use of finite loss entropy measures, as introduced in
\cite{Owen03}, allow us to go further by dropping the unnecessary
condition that the utility function has Reasonable Asymptotic
Elasticity at $+\infty$.

\section{Assumptions on $U$}

The following assumption holds throughout the paper:
\begin{assumption}\label{assumptions_on_u}
We assume that the investor has a critical wealth
$a\in[-\infty,\infty)$ and a utility function
$U:(a,\infty)\rightarrow\R$ which is increasing, strictly concave,
continuously differentiable, and satisfies the Inada conditions
\begin{equation}\label{eqn:inada}
\lim_{x\downarrow a}U'(x)=\infty,\qquad
\lim_{x\uparrow\infty}U'(x)=0.
\end{equation}
Furthermore, if the domain of $U$ is the whole real line (i.e.
$a=-\infty$) then we assume that $U$ has Reasonable Asymptotic
Elasticity at $-\infty$, in the sense that
  \[ AE_{-\infty}(U):=\liminf_{x\rightarrow-\infty}\frac{xU'(x)}{U(x)}>1. \]
This condition was introduced and discussed in detail in
\cite{Scha00}.
\end{assumption}
The convex conjugate $V$ of the utility function $U$ is defined
for $y>0$ by
\[ V(y)=\sup_{x\in(a,\infty)}\{U(x)-xy\}. \]
It is well known (cf. \cite[\S 26]{Rock72}) that under the
conditions of Assumption~\ref{assumptions_on_u}, $V$ is strictly
convex, continuously differentiable and satisfies
\begin{equation}\label{eqn:blahblah}
V'(0):=\lim_{y\downarrow0}V'(y)=-\infty\qquad\text{and}\qquad
V'(\infty):=\lim_{y\uparrow\infty}V'(y)=-a.
\end{equation}

The following Lemma, which is a simple consequence of \cite[Prop.
4.1(iii)]{Scha00}, provides an equivalent formulation of
Reasonable Asymptotic Elasticity at $-\infty$.
\begin{lemma}\label{thm:growth}
Let $U$ be a utility function defined on the whole real line, and
suppose that $U$ has reasonable asymptotic elasticity at
$-\infty$. Then there exists $b>0$ such that $V$ is positive and
increasing on $(b,\infty)$ and for any $\alpha>1$ there exists
$D>0$ such that $V(\alpha y)\le DV(y)$, for all $y\in (b,\infty)$.
\end{lemma}

\begin{proof}
Since $U$ has reasonable asymptotic elasticity at $-\infty$, a
repeated application of \cite[Prop. 4.1, (iii)]{Scha00} implies
that there exist constants $y_0>0$, $\lambda>1$ and $C>0$ such
that $V(\lambda^ny)\le C^nV(y)$, for $y\ge y_0$ and $n\in\N$.

From \eqref{eqn:blahblah} we see that $V'(\infty)=\infty$, so
there exists a $b_0>0$ such that $V$ is positive and increasing on
$(b_0,\infty)$. Set $b:=\max\{y_0,b_0\}>0$. Given any $\alpha>1$
there exists $n\in\N$ such that $\alpha\le\lambda^n$. For $y\ge b$
we have $V(\alpha y) \le V(\lambda^n y) \le C^n V(y) \le D V(y)$,
where $D:=C^n$.

\qed
\end{proof}

\section{Finite Loss Entropy Measures}

Relative to the cone $K$, we define also the convex set of {\it
separating}, or {\it pricing} measures by
\begin{equation}
M_1 := \{\Q\ll\P: X\in L^1(\Q) \textup{ and } \E[\Q]{X} \le 0
\textup{ for all }X\in K\}.
\end{equation}
In what follows, we refer frequently to the function
$V^+:=\max\{V,0\}$. Note however that in most places we can drop
the $^+$, since $V$ is convex, and its graph can be bounded from
below by a straight line.
\begin{definition} (cf. \cite{Owen03})
A measure $\Q\ll\P$ is said to have finite loss-entropy if there
exists a constant $b>0$ such that
\begin{equation}\label{eqn:lossdefn}
 \E{ V^+\left(\radnik\right)\ind_{\{\radnik\ge b\}}}<\infty.
\end{equation}
The set of pricing measures with finite loss-entropy is denoted by
\begin{equation}\label{eqn:measfinitele}
  \hatMV := \left\{\Q\in M_1: \Q\text{ has finite loss entropy}\right\}.
\end{equation}
\end{definition}

\begin{remarks}\label{remarks}
\begin{enumerate}
\item We use the~$\widehat{\phantom{.}.\phantom{.}}$~notation to
distinguish the set of finite loss entropy measures from the set
$M_V$ of finite entropy measures used in related papers. \item It
is easy to see that if \eqref{eqn:lossdefn} holds for some
constant $b>0$ then it holds for all $b>0$. In other words, if
$\Q\in\hatMV$, then $\E{V^+\left(\radnik\right)\ind_{\{\radnik\ge
b\}}}<\infty$ for all $b>0$. \item If the domain of $U$ is a
half-line, we have $\hatMV=M_1$ (see \cite[Remark 6.3]{Owen03}).
\end{enumerate}
\end{remarks}

\begin{lemma}
The set $\hatMV$ is convex.
\end{lemma}

\begin{proof}
For the case where $U$ is defined on a half real line (i.e.
$a\in(-\infty,\infty)$) convexity is trivial, since $\hatMV=M_1$.
Nevertheless, we give a universal proof. Since the function $V^+$
is convex and non-negative, given arbitrary constants $0<y\le z$,
and $x,b>0$ we have
\begin{align}\label{eqn:a1.5}
  V^+(y)\ind_{[b,\infty)}(y)
    &\le V^+(b)\ind_{[b,\infty)}(y)+V^+(z)\ind_{[b,\infty)}(y)\notag\\
  &\le V^+(x)\ind_{[b,\infty)}(x)+V^+(z)\ind_{[b,\infty)}(z)+V^+(b).
\end{align}

Take $\alpha\in[0,1]$, $\Q_1,\Q_2\in\hatMV$ and define
$\Q_\alpha:=\alpha\Q_1+(1-\alpha)\Q_2$. Let $b>0$ be an arbitrary
constant. Applying \eqref{eqn:a1.5}, we have
\begin{align*}
&\E{V^+\left(\radnik[\Q_\alpha]\right)\ind_{\{\radnik[\Q_\alpha]\ge b\}}} \\
&\qquad \le
\E{V^+\left(\radnik[\Q_1]\right)\ind_{\{\radnik[\Q_1]\ge b\}}} +
\E{V^+\left(\radnik[\Q_2]\right)\ind_{\{\radnik[\Q_2]\ge b\}}} + V^+(b)\\
& \qquad < \infty.\tag*{\qed}
\end{align*}
\end{proof}

\section{The Super-Replicable Contingent Claims}\label{sec:claims}

Let $L_U := \bigcap_{\Q\in\hatMV}L^1(\Q)$ denote the vector space
of all $\hatMV$-integrable contingent claims. Note that due to the
definition of $M_1$ it follows that $K\subseteq L_U$. Consider the
solid convex cone
\begin{equation*} K_U:=\left\{X\in L_U: X\le\widetilde X\textup{ for
some }\widetilde X\in K\right\}, \end{equation*} of all
$\hatMV$-integrable contingent claims that can be dominated by a
terminal wealth in $K$. We shall adopt throughout the common
practise of identifying probability measures $\Q\ll\P$ with their
Radon-Nikodym derivative $\radnik\in L^1_+(\P)$. Following this
convention we let $L_V^+$ denote the subspace of $L^1(\P)$, formed
by taking the linear span of the Radon-Nikodym derivatives of all
finite loss entropy pricing measures. When endowed with the
bilinear form,
\[ \ip{X}{X^+}=\E{XX^+}, \] it is easy
to show that the pair $(L_U,L_V^+)$ becomes a left dual system
(cf. \cite{Owen03}, \cite{OertOwen05}). We now define
\[
C_U:=\overline{K_U}^{\sigma(L_U,L_V^+)}.
\]
The set $C_U$ is the smallest solid, closed, convex cone in $L_U$
containing $K$, and is useful for a duality theory. We use the
subscript $U$ to stress the weak dependence of this set upon $U$.
In Section \ref{sec:representCU} we give a characterisation of
$C_U$ in terms of intersections of all $L^1(\Q)$ closures of
$K_U$. Following the terminology introduced in \cite{Owen03}, we
call $\suprep(X; C_U)$ the {\it utility-based super-replication
price} of $X$.

\section{Main Results}
Our main results are Theorems \ref{thm:duality} and
\ref{thm:main}. Recall that if $\dom(U)=\R$ then we assume only
that $U$ has Reasonable Asymptotic Elasticity at $-\infty$.
Interestingly, we need neither the assumption that $U$ is bounded,
nor that $U$ has Reasonable Asymptotic Elasticity at $+\infty$.

Given a non-empty set $A\subseteq L_U$ we let $\cone(A)$ denote
the smallest convex cone containing $A$ and we define its {\it
polar cone} $A^\lhd\subseteq L_V^+$ by
\begin{equation*}
A^\lhd:=\{X^+\in L_V^+:\ip{X}{X^+}\le0\text{ for all }X\in A\}.
\end{equation*}
For a non-empty set $B\subseteq L_V^+$ we define in a similar way
$\cone(B)$, and the polar cone $B^\lhd\subseteq L_U$ by
\begin{equation*}
B^\lhd:=\{X \in L_U:\ip{X}{X^+}\le0\text{ for all }X^+\in B\}.
\end{equation*}

We recall now, without proof, some useful facts about polar cones
(cf. \cite{OertOwen05}):
\begin{lemma}\label{thm:coneresults}
Let $A\subseteq L_U$ be non-empty. Then $A^\lhd\subseteq L_V^+$ is
a $\sigma(L_V^+,L_U)$-closed convex cone. Moreover,
\begin{equation*}
  (\cone(A))^\lhd=\overline{A}^\lhd=A^\lhd
  \qquad\text{and}\qquad
  A^{\lhd\lhd}=\overline{\cone(A)},
\end{equation*}
where closures are taken in the $\sigma(L_U,L_V^+)$ topology.
\end{lemma}

\begin{theorem}\label{thm:duality}
Suppose that $\hatMV\neq\emptyset$.  Then
\begin{equation}\label{eqn:oneway}
L_V^+ \cap \cone(M_1)=\cone(\hatMV)=K_U^\lhd=C_U^\lhd
\end{equation}
and
\begin{equation}\label{eqn:theotherway}
 C_U=(\hatMV)^\lhd.
\end{equation}
\end{theorem}
\begin{proof}
To obtain \eqref{eqn:oneway}, we first show that $L_V^+\cap M_1
\subseteq \hatMV$. To this end, take any $\Q\in L_V^+\cap M_1$.
Since $\Q$ is a probability measure, it follows from the
definition of $L_V^+$ that $\Q=\alpha\Q_0-(\alpha - 1)\Q_1$ for
some $\Q_0,\Q_1\in\hatMV$ and some $\alpha\ge1$. We now show that
$\Q\in\hatMV$:

In the case where $U$ is defined on the whole real line, let $b>0$
be the constant from the statement of Lemma~\ref{thm:growth}. Due
to Lemma~\ref{thm:growth}, we have
\begin{align*}
  \E{V^+\left(\radnik\right)\ind_{\{\radnik\ge\alpha b\}}} & \le \E{V^+\left(\alpha\radnik[\Q_0]
  \right)\ind_{\{\radnik[\Q_0]\ge b\}}} \\
  &\le D\,\E{V^+\left(\radnik[\Q_0]\right)\ind_{\{\radnik[\Q_0]\ge b\}}} <
  \infty,
\end{align*}
so $\Q\in\hatMV$. The case where $U$ is defined on a half real
line (i.e. $a\in(-\infty,\infty)$) is trivial due to Remark
\ref{remarks}(iii).

Consequently,
\[
L_V^+\cap\cone(M_1) \subseteq \cone(L_V^+\cap M_1) \subseteq
\cone(\hatMV).
\]
Now take any $X\in K_U$ and any $\Q\in\hatMV$. There exists
$\widetilde X\in K\subseteq L^1(\Q)$ such that $X\le\widetilde X$.
Hence $\mathbb E_{\Q}\big[X\big]\le\mathbb E_{\Q}\big[\widetilde
X\big]\le0$, and therefore, $\cone(\hatMV)\subseteq K_U^\lhd$.

On the other hand, since $-L_+^\infty(\P)\cup K\subseteq K_U$,
\begin{align*}
K_U^\lhd &\subseteq (-L_+^\infty(\P))^\lhd\cap K^\lhd \\
 &= \{X^+\in L^1(\P):\E{XX^+}\ge0\text{ for all }X\in L_+^\infty(\P)
 \}\cap K^\lhd\\
 &= L_+^1(\P)\cap K^\lhd \\
 &=\{X^+\in L_+^1(\P): X^+\in L_V^+\text{ and }\E{XX^+}\le0\text{ for all
 }X\in K\}\\
 &= L_V^+\cap\cone(M_1).
\end{align*}
Moreover, applying Lemma~\ref{thm:coneresults}
\begin{equation*}
C_U^\lhd=(\overline{K_U}^{\sigma(L_U,L_V^+)})^\lhd=K_U^\lhd,
\end{equation*}
and \eqref{eqn:oneway} follows. A final application of Lemma
\ref{thm:coneresults} shows that
\begin{equation*}
  C_U=\overline{K_U}^{\sigma(L_U,L_V^+)}={K_U}^{\lhd\lhd}
  =(\cone(\hatMV))^\lhd=(\hatMV)^\lhd.\tag*{\qed}
\end{equation*}
\end{proof}

\begin{theorem}\label{thm:main}
Suppose that $\hatMV\neq\emptyset$. Then for any $X\in L_U$,
\begin{equation*}
  \suprep(X;C_U) = \sup_{\Q\in \hatMV}\E[\Q]X.
\end{equation*}
\end{theorem}

\begin{proof}
Due to Theorem~\ref{thm:duality}, $C_U = (\hatMV)^\lhd$. Since
$C_U$ is solid in $L_U$ we have
\begin{align*}
\suprep(X;C_U) &=\inf\{x \in \R: X - x \in C_U \}\\
 &= \inf\{x
\in \R : \E[\Q]{X-x} \le 0
\textrm{ for all } \Q \in\hatMV\}\\
&= \inf\{x \in \R : \E[\Q]{X} \le x
\textrm{ for all } \Q \in\hatMV\}\\
&= \sup_{\Q\in\hatMV}\E[\Q]X.\tag*{\qed}
\end{align*}
\end{proof}

\begin{remark}
Given that $C_U$ is a $\sigma(L_U,L_V^+)$-closed cone, it is
specified by its polar set, $\cone(\hatMV)$. This set depends only
on the shape of $V(y)$ for arbitrarily large $y$, which in turn
depends only on the values of $U(x)$ for arbitrarily large
negative $x$. Consequently, the cone $C_U$ of allowable terminal
wealths only depends on the preferences of the investor to
asymptotically large losses. This interesting observation also
suggests an open problem: Can the set $C_U$ be parametrised by a
real number which is defined in terms of the asymptotic behaviour
of $U$ at $-\infty$?
\end{remark}

\section{A Representation of $C_U$}\label{sec:representCU}

Note that the set $K_U$ of section~\ref{sec:claims} can be
rewritten as
\begin{equation}\label{eqn:altKU}
  K_U=\bigcap\limits_{\Q\in\hatMV}\big(K - L^1_+(\Q)\big).
\end{equation}
The next theorem gives two useful alternative representations of
the weak closed cone $C_U$.
\begin{theorem}\label{thm:representCU}
\[ C_U\overset{(i)}=\bigcap_{\Q\in\hatMV}\normcl{K_U}\overset{(ii)}=\bigcap_{\Q\in\hatMV}\normcl{K-L_+^1(\Q)}. \]
\end{theorem}

\begin{proof}
(i) To show one inclusion, let
$X\in\bigcap_{\Q\in\hatMV}\normcl{K_U}$. Then for each
$\Q\in\hatMV$ there exists a sequence $\{X_n^\Q\}\subseteq K_U$
such that $X_n^\Q\overset{L^1(\Q)}\longrightarrow X$ as
$n\rightarrow\infty$. Since $K_U\subseteq(\hatMV)^\lhd$ it follows
that $\E[\Q]X=\lim_{n\rightarrow\infty}\E[\Q]{X_n^\Q}\le0$ for
each $\Q\in\hatMV$. Consequently, $X\in(\hatMV)^\lhd=C_U$.

For the other inclusion, we proceed along the lines of the proof
of the Kreps-Yan Theorem (cf. \cite[Theorem 3.5.8]{ElliKopp05})
and consider an arbitrary $Z\in L_U$ such that
$Z\not\in\normcl[\Q^*]{K_U}$ for some $\Q^*\in\hatMV$. By the Hahn
Banach Hyperplane Separation Theorem there exists a continuous
linear functional on $L^1(\Q^*)$ that separates $Z$ from the
closed cone $\normcl[\Q^*]{K_U}$. In other words there exists a
$\Lambda\in L^\infty(\Q^*)$ such that
\begin{equation}\label{eqn:aboveineq} \E[\Q^*]{\Lambda
X}\le0<\E[\Q^*]{\Lambda Z} \end{equation} for all $X\in K_U$. By
considering $X=-\ind_{\{\Lambda<0\}}\in -L_+^\infty(\P)\subseteq
K_U$, we see that $\Lambda\ge0$ $\Q^*$-a.s. and
$\E[\Q^*]\Lambda>0$. Thus, if we set
$\Lambda^*=\Lambda/\E[\Q^*]\Lambda$ then
$\Q_0(A):=\E[\Q^*]{\Lambda^*\ind_A}$ defines a probability measure
on $(\Omega,\F)$, and \eqref{eqn:aboveineq} implies that $\Q_0\in
M_1$ and $\E[\Q_0]Z>0$. To finish the proof of the first
inequality it suffices to prove that $\Q_0$ has finite loss
entropy, as then it follows from \eqref{eqn:aboveineq} that
$\Q_0\in\hatMV$ and $Z\not\in(\hatMV)^\lhd=C_U$.

In the case where $U$ is defined on the whole real line, let $b>0$
be the constant from the statement of Lemma~\ref{thm:growth}. Due
Lemma~\ref{thm:growth} it follows that
\begin{align*}
  \E{V^+\left(\radnik[\Q_0]\right)\ind_{\{\radnik[\Q_0]\ge b\|\Lambda^*\|_{L^\infty(\Q^*)}\}}}
    &=\E{V^+\left(\Lambda^*\radnik[\Q^*]\right)\ind_{\{\Lambda^*\radnik[\Q^*]\ge b\|\Lambda^*\|_{L^\infty(\Q^*)}\}}}\\
  &\le\E{V^+\left(\|\Lambda^*\|_{L^\infty(\Q^*)}\radnik[\Q^*]\right)\ind_{\{\radnik[\Q^*]\ge b\}}}\\
  &\le D\,\E{V^+\left(\radnik[\Q^*]\right)\ind_{\{\radnik[\Q^*]\ge b\}}}\\
  &<\infty.
\end{align*}
The case where $U$ is defined on a half real line (i.e.
$a\in(0,\infty)$) is trivial due to Remark \ref{remarks}(iii).

\vspace{0.5\baselineskip} \noindent(ii) To prove the second
equality it suffices to show that
\[ \normcl{K_U}=\normcl{K-L_+^1(\Q)}, \]
for an arbitrary $\Q\in\hatMV$. Indeed, from \eqref{eqn:altKU} we
have $K_U\subseteq K-L_+^1(\Q)\subseteq L^1(\Q)$, so
$\normcl{K_U}\subseteq\normcl{K-L_+^1(\Q)}$. Moreover, since
$K\cup(-L_+^\infty(\Q))\subseteq K_U$, we have
$K-L_+^\infty(\Q)\subseteq K_U$. Since $L^\infty(\Q)$ is dense in
$L^1(\Q)$, it follows that
\begin{equation*}
\normcl{K-L_+^1(\Q)}=\normcl{K-\normcl{L_+^\infty(\Q)}} \subseteq
\normcl{\normcl{K-L_+^\infty(\Q)}}\subseteq
\normcl{K_U}.\tag*{\qed}
\end{equation*}
\end{proof}

\end{document}